\documentclass[10pt, a4paper]{amsart}
\usepackage{mathrsfs}
\usepackage{amssymb, amsmath}
\usepackage{enumitem}
\usepackage{xspace}
\usepackage{pstricks,graphicx}

\newcommand{\NN}{\mathbf{N}}

\newcommand{\RR}{\mathbf{R}}

\newcommand{\cat}{{\upshape CAT($0$)}\xspace}
\newcommand{\catm}{{\upshape CAT($-1$)}\xspace}
\newcommand{\se}{\subseteq}
\newcommand{\epsi}{\epsilon}
\newcommand{\teta}{\vartheta}
\theoremstyle{plain}

\newtheorem*{thm*}{Theorem}

\theoremstyle{definition}
\newtheorem*{defn*}{Definition}

\newtheorem*{rems*}{Remarks}
\begin{document}

\title{Extreme points in non-positive curvature}

\begin{abstract}
A natural analogue of the Krein--Milman theorem is shown to fail for \cat spaces.
\end{abstract}
\author{Nicolas Monod}
\address{EPFL, 1015 Lausanne, Switzerland}
%
%
%
\maketitle
\section{Introduction}
Functional analysis sometimes offers inspiring analogies for the study of complete metric spaces of non-positive curvature such as \cat spaces. Certain classical results from functional analysis translate into fundamental facts, others into open problems, and others yet lead to counter-examples.

\medskip

Let us illustrate the first outcome on a bounded closed convex subset $C$ in a complete \cat space $X$. Typically $X$ is not locally compact and $C$ not compact. In a reflexive Banach space, such a set would still be \emph{weakly compact}. Here, $C$ is (quasi-) compact for the weakened ``convex topology''~\cite[Thm.~14]{MonodCAT0}, defined as the weakest (not necessarily Hausdorff) topology in which metrically closed \emph{convex} sets are closed.

Furthermore, in a rudimentary analogue of the Ryll-Nardzewski theorem~\cite{Ryll-Nardzewski_s}, every isometry of $C$ has a fixed point~\cite[II.2.8]{Bridson-Haefliger}.

\medskip

The second outcome occurs for instance when trying to generalize Mazur's compactness theorem~\cite{Mazur30}, because it is unknown if the closed convex hull of a compact set in $X$ remains compact~\cite[6.B$_1$(f)]{Gromov91}.

\medskip

The purpose of this note is to establish the third outcome for the Krein--Milman theorem about extreme points.

\medskip

A point in a uniquely geodesic metric space is called \emph{extreme} if it does not lie in the interior of any geodesic segment. This definition is generalized beyond unique geodesics in the presence of a \emph{bicombing}, encompassing notably all linear spaces (see e.g.~\cite{Descombes-Lang} for bicombings).

The Krein--Milman theorem, originating in~\cite{Krein-Milman}, states that a non-empty convex compact subset of a locally convex space has extreme points, and indeed is the closed convex hull of the set of extreme points. For more specific classes of Banach spaces, much stronger conclusions hold, even for bounded closed convex sets $C$ that are not compact in any weakened topology. For instance, in spaces with the Radon--Nikodym property, $C$ has extreme points and is even the closed convex hull of its strongly exposed points~\cite{Phelps74}.

\medskip
Turning back to the metric situation, it follows immediately from the \cat condition that any point \emph{maximizing the distance} to some given point must be extreme. For this reason, metrically compact \cat spaces have extreme points. In fact, in the compact case the existence of a convex bicombing is enough to retain the full conclusion of the Krein--Milman theorem~\cite{Buehler_combing}.

However, there is generally no point maximizing any distance in a bounded closed convex set. An elementary example is obtained by gluing together at a single endpoint segments of length $1-1/n$ for all $n\in \NN$. This example has, of course, plenty of extreme points: it still is the convex hull of its extreme points as in Krein--Milman. The purpose of this note is to show that it ain't necessarily so:

\begin{thm*}
There exists a bounded complete \cat space $X$ without extreme points.

Moreover, one can arrange that $X$ is compact Hausdorff for the convex topology and that every finite collection of points in $X$ is contained in a finite Euclidean simplicial complex of dimension two. Alternatively, one can construct a \catm example with hyperbolic simplicial complexes of dimension two.
\end{thm*}

Our proof uses the Pythagorean identity to play off the square-summability of the harmonic series against its non-summability. This will ensure that we remain in a finite radius while the search for extreme points can go on forever. A crucial step is to establish completeness so that there is no extreme point hiding sub rosa in the completion.

\begin{rems*}
(i)~Our example is complementary to Roberts' famous non-locally convex linear counter-example~\cite{Roberts77} to Krein--Milman: since \cat spaces satisfy a strong convexity condition, it is the non-linearity that is to blame here.

(ii)~Although there exist powerful barycentric methods for measures on \cat spaces~\cite{Korevaar-Schoen,Sturm_s}, the theorem shows that they cannot afford a Choquet theory~\cite{Choquet56}.

(iii)~Of course a space as in the theorem cannot be isometrically realized in a space where Krein--Milman holds; therefore, our ``rose'' $X$ will be folded in an appropriate way. We shall nonetheless press it flat to measure angles, obtaining an infinitely winding spiral of petals.
\end{rems*}

\section{The rose}

\begin{flushright}
\begin{minipage}[t]{0.6\linewidth}\itshape\small
C'est le temps que j'ai perdu pour ma rose\ldots
\begin{flushright}
\upshape\small
Antoine de Saint-Exup\'ery, \emph{Le Petit Prince}, 1943
\end{flushright}
\end{minipage}
\end{flushright}
\vspace{3mm}

Define for each $n\in\NN$ the radius $r_n= \sqrt{\sum_{p=1}^n 1/p^2 }$. The double \emph{petal} $P_n\se \RR^2$ shall be the closed set
$$P_n = \Big\{ (x,y) : (n+1) r_n |x| \leq |y|\leq r_n \Big\}$$
with the path-metric induced from $\RR^2$. In particular, $P_n$ is a compact \cat space; it can also be viewed as two Euclidean triangles glued at a tip.

We refer to the locus $x=0$ as the \emph{central segment} of $P_n$; it has length $2 r_n$. We call the loci $|y|=(n+1) r_n x$ and $|y|= -(n+1) r_n x$ the \emph{right border} and \emph{left border} of $P_n$; they are segments of length $2r_{n+1}$.

\begin{figure*}[!h]
\centering{
\resizebox{0.45\textwidth}{!}{\input{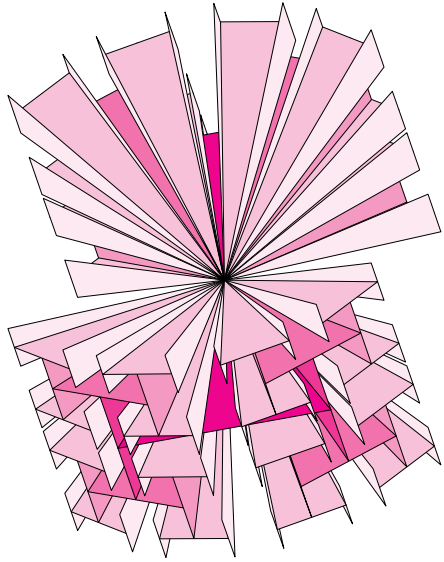}}}
\caption{The first petals, up to $P_6$, of the (Euclidean) rose.}\label{fig:3D}
\end{figure*}

The \emph{rose} is the \cat space $X$ obtained by inductively gluing $2^{n-1}$ copies of $P_n$ as follows, starting from a single copy of $P_1$. For each $n\in \NN$ and for each copy of $P_n$ we assign two new copies of $P_{n+1}$. The first is glued by identifying its central segment with the right border of the given copy of $P_n$. The second is glued by identifying its central segment with the left border. The gluings make sense since these segments have all length $2 r_{n+1}$, and $X$ is defined to be the increasing union of the successively glued spaces (Figure~\ref{fig:3D}).

This rose $X$ is indeed a \cat space, see~\cite[II.11.1]{Bridson-Haefliger}. Notice that each verification of the \cat condition takes place in a finite gluing since the convex hull of a set of point does not visit petals of higher index than these points. The rose is bounded because every point of $P_n$ is at distance at most $r_{n+1}$ of the \emph{center}, which is the common intersection of all copies of all petals. Thus $X$ has radius $\pi/\sqrt 6$.

Furthermore, by construction $X$ has no extreme point. Indeed, if a point in some $P_n$ were extreme in $X$, it would a fortiori be extreme in $P_n$ and hence it would be one of the four outer corners of $P_n$. These four points are, however, all midpoints of segments in appropriate copies of $P_{n+1}$.

\smallskip
We now turn to the critical point and prove that $X$ is complete.

Let $(x_j)_{j=1}^\infty$ be a Cauchy sequence in $X$. Suppose for a contradiction that $(x_j)$ does not converge. Upon discarding finitely many terms, we can assume that $(x_j)$ remains at distance at least $\epsi$ of the center for some $\epsi>0$. Define the \emph{index} $n(x)$ of a point $x\in X$ to be the smallest $n$ such that some copy of $P_n$ contains $x$; if $x$ is not the center, then it is contained in at most one other petal, namely a copy of $P_{n(x)+1}$. The indices $n(x_j)$ are unbounded as $j$ varies, because otherwise $(x_j)$ would be confined to the gluing of finitely many petals, which is a compact metric space; this would imply that $(x_j)$ converges.

In order to contradict the Cauchy assumption, it suffices to prove that for any $j$ there is $i>j$ with $d(x_j, x_i) \geq \epsi$. Fix thus $j$ and consider only those $i>j$ with $n(x_i) \geq n(x_j) +2$.

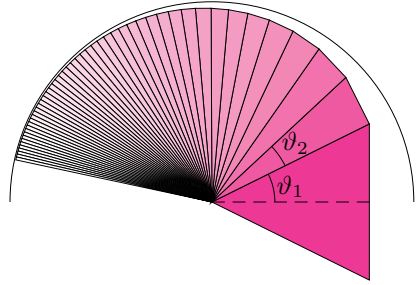
\begin{figure*}[!h]
\centering{
\resizebox{0.5\textwidth}{!}{\psset{unit=2.1cm}
\begin{pspicture}(-1.5,-0.5)(1.5,1.5)
\psarc[linewidth=0pt](0,0){1.28254983}{0}{180}
\pspolygon[linewidth=0pt,opacity=1,fillcolor=magenta!90,fillstyle=solid](1,0.5)(1,-0.5)(0,0)
\psline[linewidth=0pt,opacity=1,linestyle=dashed](0,0)(1,0)
\pspolygon[linewidth=0pt,opacity=1,fillcolor=magenta!80,fillstyle=solid](1,0.5)(0.8509288015,0.798142397)(0,0)
\pspolygon[linewidth=0pt,opacity=1,fillcolor=magenta!70,fillstyle=solid](0.8509288015,0.798142397)(0.6798982879,0.980484283)(0,0)
\pspolygon[linewidth=0pt,opacity=1,fillcolor=magenta!60,fillstyle=solid](0.6798982879,0.980484283)(0.5155463052,1.094451058)(0,0)
\pspolygon[linewidth=0pt,opacity=1,fillcolor=magenta!55,fillstyle=solid](0.5155463052,1.094451058)(0.3647702953,1.1654748048)(0,0)
\pspolygon[linewidth=0pt,opacity=1,fillcolor=magenta!50,fillstyle=solid](0.3647702953,1.1654748048)(0.2284346563,1.208145132)(0,0)
\pspolygon[linewidth=0pt,opacity=1,fillcolor=magenta!49,fillstyle=solid](0.2284346563,1.208145132)(0.1056109003,1.2313685029)(0,0)
\pspolygon[linewidth=0pt,opacity=1,fillcolor=magenta!48,fillstyle=solid](0.1056109003,1.2313685029)(-0.0050937853,1.2408633223)(0,0)
\pspolygon[linewidth=0pt,opacity=1,fillcolor=magenta!46,fillstyle=solid](-0.0050937853,1.2408633223)(-0.1050929428,1.2404528224)(0,0)
\pspolygon[linewidth=0pt,opacity=1,fillcolor=magenta!44,fillstyle=solid](-0.1050929428,1.2404528224)(-0.1956775195,1.2327783671)(0,0)
\pspolygon[linewidth=0pt,opacity=1,fillcolor=magenta!42,fillstyle=solid](-0.1956775195,1.2327783671)(-0.2779804966,1.2197145084)(0,0)
\pspolygon[linewidth=0pt,opacity=1,fillcolor=magenta!40,fillstyle=solid](-0.2779804966,1.2197145084)(-0.3529804352,1.2026215575)(0,0)
\pspolygon[linewidth=0pt,opacity=1,fillcolor=magenta!38,fillstyle=solid](-0.3529804352,1.2026215575)(-0.4215178212,1.1825052074)(0,0)
\pspolygon[linewidth=0pt,opacity=1,fillcolor=magenta!36,fillstyle=solid](-0.4215178212,1.1825052074)(-0.4843141521,1.1601207202)(0,0)
\pspolygon[linewidth=0pt,opacity=1,fillcolor=magenta!34,fillstyle=solid](-0.4843141521,1.1601207202)(-0.5419900223,1.1360428465)(0,0)
\pspolygon[linewidth=0pt,opacity=1,fillcolor=magenta!32,fillstyle=solid](-0.5419900223,1.1360428465)(-0.5950809973,1.110713891)(0,0)
\pspolygon[linewidth=0pt,opacity=1,fillcolor=magenta!30,fillstyle=solid](-0.5950809973,1.110713891)(-0.6440510845,1.08447746)(0,0)
\pspolygon[linewidth=0pt,opacity=1,fillcolor=magenta!28,fillstyle=solid](-0.6440510845,1.08447746)(-0.6893039851,1.0576026002)(0,0)
\pspolygon[linewidth=0pt,opacity=1,fillcolor=magenta!26,fillstyle=solid](-0.6893039851,1.0576026002)(-0.7311924299,1.0303013513)(0,0)
\pspolygon[linewidth=0pt,opacity=1,fillcolor=magenta!24,fillstyle=solid](-0.7311924299,1.0303013513)(-0.7700259108,1.0027416987)(0,0)
\pspolygon[linewidth=0pt,opacity=1,fillcolor=magenta!23,fillstyle=solid](-0.7700259108,1.0027416987)(-0.8060770899,0.9750572591)(0,0)
\pspolygon[linewidth=0pt,opacity=1,fillcolor=magenta!22,fillstyle=solid](-0.8060770899,0.9750572591)(-0.8395871233,0.9473546089)(0,0)
\pspolygon[linewidth=0pt,opacity=1,fillcolor=magenta!21,fillstyle=solid](-0.8395871233,0.9473546089)(-0.870770101,0.919718889)(0,0)
\pspolygon[linewidth=0pt,opacity=1,fillcolor=magenta!20,fillstyle=solid](-0.870770101,0.919718889)(-0.8998167616,0.8922181343)(0,0)
\pspolygon[linewidth=0pt,opacity=1,fillcolor=magenta!19,fillstyle=solid](-0.8998167616,0.8922181343)(-0.9268976139,0.8649066463)(0,0)
\pspolygon[linewidth=0pt,opacity=1,fillcolor=magenta!18,fillstyle=solid](-0.9268976139,0.8649066463)(-0.9521655718,0.8378276431)(0,0)
\pspolygon[linewidth=0pt,opacity=1,fillcolor=magenta!17,fillstyle=solid](-0.9521655718,0.8378276431)(-0.9757581864,0.8110153559)(0,0)
\pspolygon[linewidth=0pt,opacity=1,fillcolor=magenta!16,fillstyle=solid](-0.9757581864,0.8110153559)(-0.997799546,0.7844967001)(0,0)
\pspolygon[linewidth=0pt,opacity=1,fillcolor=magenta!15,fillstyle=solid](-0.997799546,0.7844967001)(-1.0184018997,0.7582926139)(0,0)
\pspolygon[linewidth=0pt,opacity=1,fillcolor=magenta!14,fillstyle=solid](-1.0184018997,0.7582926139)(-1.0376670496,0.7324191386)(0,0)
\pspolygon[linewidth=0pt,opacity=1,fillcolor=magenta!13,fillstyle=solid](-1.0376670496,0.7324191386)(-1.0556875502,0.7068882933)(0,0)
\pspolygon[linewidth=0pt,opacity=1,fillcolor=magenta!12,fillstyle=solid](-1.0556875502,0.7068882933)(-1.0725477447,0.6817087881)(0,0)
\pspolygon[linewidth=0pt,opacity=1,fillcolor=magenta!11,fillstyle=solid](-1.0725477447,0.6817087881)(-1.0883246643,0.6568866062)(0,0)
\pspolygon[linewidth=0pt,opacity=1,fillcolor=magenta!10,fillstyle=solid](-1.0883246643,0.6568866062)(-1.1030888101,0.6324254834)(0,0)
\pspolygon[linewidth=0pt,opacity=1,fillcolor=magenta!9,fillstyle=solid](-1.1030888101,0.6324254834)(-1.1169048364,0.6083273018)(0,0)
\pspolygon[linewidth=0pt,opacity=1,fillcolor=magenta!9,fillstyle=solid](-1.1169048364,0.6083273018)(-1.1298321496,0.5845924168)(0,0)
\pspolygon[linewidth=0pt,opacity=1,fillcolor=magenta!8,fillstyle=solid](-1.1298321496,0.5845924168)(-1.1419254328,0.5612199271)(0,0)
\pspolygon[linewidth=0pt,opacity=1,fillcolor=magenta!7,fillstyle=solid](-1.1419254328,0.5612199271)(-1.1532351102,0.5382079)(0,0)
\pspolygon[linewidth=0pt,opacity=1,fillcolor=magenta!6,fillstyle=solid](-1.1532351102,0.5382079)(-1.1638077552,0.5155535587)(0,0)
\pspolygon[linewidth=0pt,opacity=1,fillcolor=magenta!5,fillstyle=solid](-1.1638077552,0.5155535587)(-1.1736864536,0.4932534395)(0,0)
\pspolygon[linewidth=0pt,opacity=1,fillcolor=magenta!4,fillstyle=solid](-1.1736864536,0.4932534395)(-1.1829111252,0.4713035225)(0,0)
\pspolygon[linewidth=0pt,opacity=1,fillcolor=magenta!3,fillstyle=solid](-1.1829111252,0.4713035225)(-1.1915188104,0.4496993415)(0,0)
\pspolygon[linewidth=0pt,opacity=1,fillcolor=magenta!3,fillstyle=solid](-1.1915188104,0.4496993415)(-1.1995439261,0.4284360764)(0,0)
\pspolygon[linewidth=0pt,opacity=1,fillcolor=magenta!2,fillstyle=solid](-1.1995439261,0.4284360764)(-1.2070184943,0.4075086305)(0,0)
\pspolygon[linewidth=0pt,opacity=1,fillcolor=magenta!1,fillstyle=solid](-1.2070184943,0.4075086305)(-1.2139723469,0.3869116956)(0,0)
\pspolygon[linewidth=0pt,opacity=1,fillcolor=magenta!1,fillstyle=solid](-1.2139723469,0.3869116956)(-1.2204333103,0.3666398067)(0,0)
\pspolygon[linewidth=0pt,opacity=1,fillcolor=magenta!1,fillstyle=solid](-1.2204333103,0.3666398067)(-1.2264273707,0.3466873879)(0,0)
\pspolygon[linewidth=0pt,opacity=1,fillcolor=magenta!1,fillstyle=solid](-1.2264273707,0.3466873879)(-1.2319788237,0.3270487908)(0,0)
\pspolygon[linewidth=0pt,opacity=1,fillcolor=magenta!0,fillstyle=solid](-1.2319788237,0.3270487908)(-1.2371104093,0.307718327)(0,0)
\pspolygon[linewidth=0pt,opacity=1,fillcolor=magenta!0,fillstyle=solid](-1.2371104093,0.307718327)(-1.241843434,0.2886902953)(0,0)
\pspolygon[linewidth=0pt,opacity=1,fillcolor=magenta!0,fillstyle=solid](-1.241843434,0.2886902953)(-1.2461978815,0.269959004)(0,0)
\psarc[linewidth=0pt](0,0){0.4}{0}{26.5651}
\rput(0.5,0.1){{$\vartheta_1$}}
\psarc[linewidth=0pt](0,0){0.52}{26.5651}{43.1666}
\rput(0.53,0.38){{$\vartheta_2$}}
\end{pspicture}}}
\caption{The petals up to $P_{50}$, folded and pressed flat in a circle of radius $\pi/\sqrt 6$; ultimately, the $P_n$ wind around infinitely often.}\label{fig:flat}
\end{figure*}

We can assume that the segment $[x_j, x_i]$ avoids the center of the rose since otherwise $d(x_j, x_i) \geq 2 \epsi$. It follows that $[x_j, x_i]$ traverses notably a non-trivial portion of successive copies of $P_{n(x_j)+1}, P_{n(x_j)+2}, \ldots, P_{n(x_i)-1}$, each time entering on the central segment and leaving through a border. Thus, the Alexandrov angle formed at the center of $X$ by $x_j$ and $x_i$ is at least $\sum_{n=n(x_j)+1}^{n(x_i)-1} \teta_n$, where $\teta_n$ is the angle between the central segment and the borders of $P_n$ (Figure~\ref{fig:flat}). That angle satisfies $\teta_n > \sin \teta_n > \sqrt 6/\pi(n+1)$. It now follows from the divergence of the harmonic series that if we choose $i$ so that $n(x_i)$ is large enough compared to $n(x_j)$, then the angle at the center between $x_j$ and $x_i$ will be at least $\pi/3$ and hence $d(x_j, x_i) \geq \epsi$.

We observe in passing that the above argument with $\pi/3$ replaced by $\pi$ shows in fact that $[x,y]$ contains the center as soon as $n(y)$ is large enough compared to $n(x)$. Therefore, any sequence $(y_j)$ with $n(y_j)$ going to infinity will converge to the center in the convex topology.

\smallskip
To prove that the convex topology of $X$ is Hausdorff amounts to the following. Given distinct points $x,y\in X$, we need to cover $X$ with finitely many closed convex subsets $U_i, V_i\se X$ such that $x\notin U_i$ and $y\notin V_i$ for all $i$. We shall apply the following sufficient condition: it is enough to find a (metrically) compact convex subset $K$ containing $x,y$ such that $X\setminus K$ can be covered by finitely many closed convex subsets $W_i\se X$ with $y\notin W_i$ for all $i$. This is indeed sufficient because we can first cover $K$ with finitely many closed convex subsets $U_i, V_i\se K$ as above since the convex topology of $K$ coincides with the metric topology~\cite[Lem.~17]{MonodCAT0}, hence is Hausdorff. These sets $U_i, V_i$ are still closed and convex as subsets of $X$ and it now suffices to add all $W_i$ to the collection of $V_i$ to cover $X$.

We now verify the condition. Upon possibly exchanging $x$ and $y$ we can assume that $y$ is not the center. Choose $n>n(x), n(y)$ and define $K$ to be the union of all copies of $P_m$ over all $m\leq n$ with the gluings of the rose. Then define $W_1, \ldots, W_{2^n}$ to be the $2^n$ connected components of the union of all copies of $P_m$ over $m>n$, with the gluings specified by the construction for $m>n$ only. These sets are all closed convex and satisfy the criteria of the sufficient condition.

\smallskip
This concludes the proof of the theorem for the \cat statement. The \catm construction is virtually identical with hyperbolic triangles, the key being that the trigonometric estimates remain the same at the first order as the angles $\teta_n$ converge to zero.

\medskip
We observe that in either case the full isometry group of the rose is an infinite iterated wreath product of Klein four-groups (see e.g.~\cite[\S IV.4]{BOERT} for such iterated products).



\bibliographystyle{../BIB/amsalpha}
\bibliography{../BIB/ma_bib}

\end{document}